\documentclass[12pt, final]{l4dc2021} 
\makeatletter
\def\set@curr@file#1{\def\@curr@file{#1}} 
\makeatother
\usepackage[load-configurations=version-1]{siunitx} 
\title[Stable Online Control of LTV Systems]{Stable Online Control of Linear Time-Varying Systems}
\usepackage{times}
\usepackage{float}
\usepackage{wrapfig} 
\usepackage{xcolor}
\usepackage{enumitem}
\graphicspath{ {./images/} }
\SetKwInOut{Parameter}{Parameters}
\date{}

\newcommand{\revise}[1]{{\color{black} #1 }}
\newcommand*{\qed}{\hfill\ensuremath{\blacksquare}}%

\DeclareMathOperator{\Tr}{Tr}

\allowdisplaybreaks

\author{%
 \Name{Guannan Qu}$^*$ \Email{gqu@caltech.edu}\\%
 \Name{Yuanyuan Shi}$^*$ \Email{yshi7@caltech.edu}\\%
 \Name{Sahin Lale}$^*$ \Email{alale@caltech.edu}\\%
 \Name{Anima Anandkumar} \Email{anima@caltech.edu}\\%
 \Name{Adam Wierman} \Email{adamw@caltech.edu}\\%
 \addr California Institute of Technology, Pasadena, CA \\
 $^*$Equal contribution}

\begin{document}

\maketitle

\begin{abstract}%
 Linear time-varying (LTV) systems are widely used for modeling real-world dynamical systems due to their generality and simplicity. Providing stability guarantees for LTV systems is one of the central problems in control theory. However, existing approaches that guarantee stability typically lead to significantly sub-optimal cumulative control cost in online settings where only current or short-term system information is available. In this work, we propose an efficient online control algorithm, COvariance Constrained Online Linear Quadratic (COCO-LQ) control, that  guarantees input-to-state stability for a large class of LTV systems while also minimizing the control cost. The proposed method incorporates a state covariance constraint into the semi-definite programming (SDP) formulation of the LQ optimal controller. We empirically demonstrate the performance of COCO-LQ in both synthetic experiments and a power system frequency control example. 
\end{abstract}


\section{Introduction}

Time-invariant systems have traditionally been the main focus of the study for the linear dynamical systems community. 
However, real-world systems are often \emph{time-varying.} For example, consider a power system that includes renewable generation (e.g. solar/wind). Due to the intermittency of renewable energy, the system dynamics for frequency regulation in the power system are time-varying. Applying a time-invariant controller in this setting may lead to frequency instability and line failures~\citep{ulbig2014impact}. Time-varying systems are also crucial for many other applications, such as autonomous vehicles and aircraft control~\citep{falcone2008linear}. While not all time-varying systems have linear dynamics, many applications can be approximated by linear time-varying (LTV) systems via a local linear approximation at each time step~\citep{todorov2005generalized}, e.g., the frequency control example described above.  As a result, LTV systems are widely-used and there is a large literature focused on designing controllers for LTV systems~\citep{amato2010finite, ouyang2017learning}.  

Perhaps the most fundamental challenge in dynamical systems is stability.  While the design of stable linear time-invariant (LTI) systems is well understood, the same cannot be said for LTV systems. To this point, several notions of stability have recieved attention, e.g., input-to-state stability (ISS), mean-square stability and Lyapunov stability. 
ISS is the most widely adopted notion and  aims to guarantee the boundedness of the state given bounded initial conditions \citep{hong2010finite}. In most applications of LTV systems, it is crucial to guarantee ISS both in order to avoid saturation, maintain the robustness and validity of linearization \citep{tarbouriech2006advanced,khalil2002nonlinear}. 

While there is considerable prior work focused on stability in LTV systems, most prior work studies stability in the offline setting where either the sequence of system parameters are known, e.g., \citep{amato2010finite, li2019lyapunov}, or the system parameters have a particular variation pattern, e.g.,  \citep{garcia2009finite}. Maintaining ISS guarantees becomes significantly harder in the online setting where the system parameters are observed in real-time and may have arbitrary variations. This online setting is the most relevant to many applications, e.g., frequency regulation.

Though stability is crucial, it is not enough for a controller to be stable. A controller must also have low cost.  For instance, in order to stabilize the dynamics, a controller may use arbitrarily big control inputs, which may result in sub-optimal cost.   
In classical optimal control problems, e.g. the time-varying linear quadratic (LQ) control setting, the goal is to design a stabilizing controller that minimizes the cost for a particular finite horizon while assuming access to the whole trajectory for that duration.  It is possible to characterize the optimal policy in such settings~\citep{bertsekas1995dynamic}; however, in the online setting when only current or short-termed system information is available, these methods
may not guarantee stability, e.g., see Section \ref{sec:naive}. There have been recent efforts to provide sub-optimality guarantees on the acquired cost in the online LTV setting, e.g., \citep{gradu2020adaptive}, but it is unclear if the proposed controllers maintain stability for all time-steps since the main focus is on minimizing the cumulative cost.

Thus, despite considerable recent work, much remains to be understood about the design of online LTV controllers.  In particular, this paper is motivated by the following question: 
\begin{center}
    \emph{Is it possible for an online controller to guarantee stability and maintain low cost in LTV systems?}
\end{center}

\paragraph{Contributions.}

In this work, we answer question above affirmatively. Specifically, we propose \textbf{Co}variance \textbf{C}onstrained \textbf{O}nline \textbf{L}inear \textbf{Q}uadratic (COCO-LQ) control, a novel online control algorithm that aims to minimize the control cost while ensuring provable stability guarantees in LTV systems without restricting how slow or fast the underlying system changes. Further, we demonstrate the performance of the proposed method in various synthetic LTV systems and in the power system frequency control example that motivated our study.

The main technical contribution of the paper is a stability guarantee for COCO-LQ in LTV systems. Specifically, we show that COCO-LQ guarantees ISS in online time-varying systems. The key technique that underpins the proposed algorithm is the addition of a novel semi-definiteness constraint on the state covariance matrix into the standard online semi-definite programming (SDP) formulation of linear quadratic optimal control. We show that this constraint promotes the sequential strong stability of the controllers~\citep{cohen2018online}, which in turn guarantees ISS with a proper choice of an algorithm hyperparameter. Adding this additional constraint is simple and does not result in a significant increase of computational complexity compared to the standard LQ formulation. Moreover, we prove that if the proposed SDP is not directly feasible,  short-term predictions on the future system parameters are necessary and can be used in COCO-LQ in order to ensure ISS. 

\paragraph{Related work.} The work in this paper builds on the design of linear time-invariant (LTI) controllers to provide a new approach for the design of stable controllers for linear-time-varying (LTV) systems.  As such, we describe related work on both LTI and LTV systems below.

\textit{LTI Systems.} In study of control of LTI systems, linear quadratic regulator (LQR) 
has been considered in detail. In the classical setting where the underlying system is known, the optimal control law is given by a linear feedback controller obtained by solving Riccati equations \citep{bertsekas1995dynamic}.  Alternatively, the optimal control problem can also be posed via semi-definite programming (SDP) \citep{vandenberghe1996semidefinite}, which is the approach we build on in the current paper. 

Recently, there has been growing interest in online control of these linear systems when the underlying dynamics are unknown. Most of these works study the problem with a regret minimization perspective, e.g.,~\citep{abbasi2011regret,dean2018regret,lale2020adaptive,lale2020explore}.
However, these methods have so far only been applied in LTI systems with time-varying costs and disturbances. Extensions to LTV dynamics, which are the focus of this paper, are not known. 

\textit{LTV Systems.} As in the case of LTI systems, optimal control of LTV systems where the sequence of system parameters can be obtained by solving backwards Riccati equations \citep{bertsekas1995dynamic}. However, in the online case when the sequence of systems is unknown, the design of controllers is challenging. There are several lines of work in adaptive control and model-predictive control (MPC) that have been studied to this point. In adaptive control of LTV systems, the underlying systems are unknown and the results generally assume slow and bounded or fixed systematic variation of dynamics with bounded disturbances \citep{middleton1988adaptive,marino2000robust,ouyang2017learning}. In MPC of LTV systems, a finite horizon of sequence of systems (predictions) is known and the system is again assumed to be slowly varying or open-loop stable, e.g., \citep{zheng1994robust, falcone2007linear}.
Different from prior works, in the current work we consider the online problem and make no assumptions about how the system varies over time.

As in the LTI setting, the study of regret minimization in LTV systems has recently received attention.~\citet{goel2020regret,gradu2020adaptive} are most related to the current paper. \citet{goel2020regret} considers the setting where the sequence of systems is known 
and provides regret-optimal controller framework. \citet{gradu2020adaptive} studies the adaptive regret of online control in LTV systems with bounded cost. Note that when the cost is bounded, a finite regret need not guarantee stability. In contrast, we use a quadratic (unbounded) cost and we can guarantee stability. 



\paragraph{Notation.} We denote the Euclidean norm of a vector $x$ as $\|x\|$. For a matrix $A$, $\| A \|$ is its spectral norm, $A^\top$ is its transpose, and $\Tr(A)$ is its trace. 
$\mathcal{N}(\mu, \Sigma)$ denotes normal distribution with mean $\mu$ and covariance $\Sigma$. $A \succ B$ and $A \succeq B$ denote that $A-B$ is positive definite and positive semi-definite respectively. $A \bullet B$ denotes the element-wise inner product of $A$ and $B$, \textit{i.e.}, $\Tr(A^\top B)$. 

\section{Model \& Background}



We consider the following linear time-varying (LTV) system, 
\begin{equation}\label{eq:lti_dyn}
    x_{t+1} = A_t x_t + B_tu_t + w_t,
\end{equation}
where $x_t \in \mathbb{R}^{d}$ is the system state, $u_t \in \mathbb{R}^{p}$ is the control input and $w_t \in \mathbb{R}^{d}$ is the disturbance at time $t$. The system is stochastic, \textit{i.e.}, $w_t \sim \mathcal{N}(0, W)$ for $W\succ 0$. The cost at each time-step is a quadratic function of the state and control, $x_t^{\top}Q x_t+u_t^{\top}R u_t$, where $Q, R \succ 0$. 

The decision maker operates in an online setting. That is, at each time-step $t$, the learner observes the state $x_t$ and system matrix $(A_t, B_t)$ before choosing action $u_t$  and suffering cost $x_t^{\top}Qx_t+u_t^{\top}Ru_t$. We assume that the cost matrices $(Q, R)$ are time-invariant and known to the learner.  However, future system matrices $(A_{t+1}, ..., A_T)$ and $(B_{t+1}, ..., B_T)$ are unknown to the learner and are chosen by the environment, potentially stochastically or adversarially.



\paragraph{Stability.} One of the most central goals for controller design is to ensure stability. In this work, we focus on the notion of input to state stability (ISS) and strive to design controllers that provide ISS. ISS has been the main notion of stability considered in designing stabilizing controllers both in linear and nonlinear systems ~\citep{hong2010finite,sontag2008input,jiang2001input}. To formally define ISS, let $\mathcal{K}_\infty$ be the set of functions from nonnegative reals to nonnegative reals that are continuous, strictly increasing, bijective. Then, ISS is defined as follows.
\begin{definition}[ISS]
A LTV system with deterministic policy $\mathcal{A}$ is said to be input to state stable if there exists functions $\beta_1:[0,\infty)\times \mathbb{N}\rightarrow [0,\infty)$ and $\beta_2\in\mathcal{K}_\infty$ that satisfy $\beta_1(\cdot ,t) \in\mathcal{K}_\infty$ for any $t\in \mathbb{N}$, $\lim_{t\rightarrow \infty} \beta_1(a,t) = 0$ for any $a\geq 0$ such that, for any disturbance sequence $\{w_t\}_{t=0}^\infty$, any initial time $t_0$, any initial state $x_{t_0}$, and any $t\geq t_0$, we have
$    \Vert x_t\Vert \leq \beta_1(\Vert x_{t_0}\Vert ,t-t_0) + \beta_2(\sup_{t'\in \mathbb{N}} \Vert w_{t'}\Vert)$.
\end{definition}


\paragraph{Cost.} In addition to stability, another important objective for controller design is  maintaining a small, near-optimal control cost. Here we adopt the standard linear quadratic (LQ) cost model, \textit{i.e.},  
\revise{\begin{align}
    J_{T}(\mathcal{A}) = \lim_{T \rightarrow \infty} \frac{1}{T} \mathbb{E} \left[\sum\nolimits_{t=1}^{T} x_t^{\top}Qx_t+u_t^{\top} R u_t\right], \label{eq:lq_cost}
\end{align}}
\noindent where $u_1, \ldots, u_t$ are chosen according to policy $\mathcal{A}$, and the expectation is taken with respect to randomness of noise sequence $w_t$.

In this work, our goal is to ensure both stability and near-optimal cost. It should be noted that there is a trade-off between these two goals. On the one hand, a stabilizing controller without cost-awareness may produce arbitrarily large control inputs and induce high cost, which is impractical to implement. On the other hand, a greedy approach that merely focuses on cost minimization may lead to instability, as we highlight in the Section~\ref{sec:naive} below.


Though our focus is on LTV systems, our approach builds on the SDP formulation of the optimal controller for LTI systems in \citep{vandenberghe1996semidefinite}.
\begin{proposition}\citep{vandenberghe1996semidefinite}\label{prop:optimal_lqr_sdp}
When $A_t = A, B_t = B$ and $(A,B)$ is controllable, the optimal $K^* = LQR(A,B,Q,R)$ where $u_t=K^*x_t$, can be obtained by the following SDP
\begin{align*}
        \min_{\Sigma \succeq 0} &\begin{bmatrix}Q & 0\\0 &R
    \end{bmatrix} \bullet\Sigma\qquad
    \text{s.t.}\qquad \Sigma_{xx} = \begin{bmatrix} A_t & B_t\end{bmatrix} \Sigma \begin{bmatrix} A_t & B_t\end{bmatrix}^{\top}+W,
\end{align*}
which has a unique symmetric solution $\Sigma^*$ that decomposes to the following blocks
$ \Sigma^* = \begin{bmatrix}
\Sigma_{xx}^* & \Sigma_{xu}^* \\ {\Sigma_{xu}^*}^{\top} & \Sigma_{uu}^*  \end{bmatrix}, $
where $\Sigma_{xx}^* \!\!\in\! \mathbb{R}^{d \times d}$, $\Sigma_{xu}^* \!\!\in\!  \mathbb{R}^{d \times p}$ and $\Sigma_{uu}^* \!\!\in\! \mathbb{R}^{p \times p}$. Then, the optimal controller is $K^* \!=\! {\Sigma_{xu}^*}^{\top} (\Sigma_{xx}^*)^{-1}$. 
\end{proposition}

\noindent The optimal LQR controller described above both stabilizes the system and achieves the minimum cost.  The current paper makes a step toward understanding if it is possible to extend this formulation to the case of LTV systems.




\section{A Naive Approach}\label{sec:naive}
How to achieve stable, cost-optimal control of LTI systems is well-known; however this is not the case in LTV systems. To illustrate the challenge of online control of LTV systems, we start by studying the performance of 
a naive ``plug in'' approach where upon receiving $(A_t, B_t)$ an optimal controller for $A_t, B_t$ is computed under the assumption that the system is time-invariant. Due to its simplicity, this approach has been employed in many contexts, e.g. \citet{li2019online_mdp} for a MDP setting.  
In this section we provide an example which shows that such a myopic approach based on optimal LTI control described above fails to stabilize the system even in simple settings where $A_t$ can only switch between two possible choices and $B_t$ is fixed. This highlights that one cannot naively apply LTI design approaches in LTV systems and expect to maintain stability.  

\begin{example}\label{example}
Consider a system with $Q =\epsilon I$, $R =  I$, $w_t = 0$, and
$A = \left[\begin{array}{cc}
     \rho &0  \\
        a  &\rho 
\end{array} \right], A' = \left[\begin{array}{cc}
     \rho & a  \\
        0  &\rho 
\end{array} \right] ,$
where $0\!<\!\rho\!<\!1$, and $a \!>\!\sqrt{2}$. Suppose $A_t$ alternates between $A$ and $A'$ and $B_t\!=\! B\!=\! I$. Define the optimal LTI controllers for $A$ and $A'$ as
$ K \!:=\! LQR(A,B,Q,R)$ and $K' \!:=\! LQR(A',B,Q,R). $

To show that the optimal LTI controllers will not stabilize the system, we consider a case where $\epsilon\to 0$. In this case, one can check that $K,K' \rightarrow 0$, and more precisely, $K,K' = O(\epsilon)$. 
Since $A_t$ alternates between $A,A'$,  $K_t$ also alternates between $K$ and $K'$ under the myopic design we are considering.  Thus, the system state follows 
$x_{t+2} = (A+K)(A'+K') x_t.$ Notice that as $\epsilon \rightarrow 0$,
$(A + K)(A'+K') \rightarrow A A' = \left[\begin{array}{cc}
     \rho^2 & a\rho  \\
        a\rho  & a^2 + \rho^ 2
\end{array} \right]$. 
Here, $AA'$ is unstable since its largest eigenvalue is greater than $\frac{1}{2}\Tr(AA') \!=\! \rho^2 \!+\! \frac{a^2}{2} >1$. Thus, for small enough $\epsilon$, the naive strategy that uses the LTI controller at each time-step leads to instability.

\end{example}

\section{Main Result} 
The previous section highlights that a naive application of LTI control cannot guarantee stability for LTV systems. We now propose a new approach, COvariance Constrained Online LQ (COCO-LQ) control (Section~\ref{subsec:algorithm}).
Our main technical result shows that COCO-LQ provably guarantees stability in LTV systems (Section~\ref{subsec:guarantee}) when the SDP is feasible. In Section~\ref{subsec:infeasibility}, we discuss how to handle the situation when the SDP is infeasible and Section~\ref{sec:estimation_error} discuss the effect of model estimation error. Detailed proofs could be found in the Appendix.
\vspace{-6pt}
\subsection{COvariance Constrained Online LQ (COCO-LQ)}\label{subsec:algorithm}
The naive approach discussed in Section~\ref{sec:naive} seeks to solve the LTI problem at every time step, which is equivalent to solving the SDP in Proposition~\ref{prop:optimal_lqr_sdp} for every $(A_t, B_t)$. The reason this method fails is that it only considers cost minimization without explicitly considering stability. 
The main idea of COCO-LQ is to enforce stability via a state covariance constraint embedded into the SDP framework. The proposed algorithm is stated formally in Algorithm \ref{alg1}.
\setlength{\textfloatsep}{6pt}
\begin{algorithm2e}[t]
	\SetAlgoLined
	\DontPrintSemicolon
	\LinesNumbered
	\Parameter{$\alpha \in [0, 1)$}
	\KwIn{$Q, R, W \succ 0$}
	\For{$t=1, 2, ...$}{
		\textbf{Receive} state $x_t$, and system parameter $A_t, B_t$\;
		\textbf{Compute policy:} Let $\Sigma_{t} \in \mathbb{R}^{n \times n}$ be an optimal solution to the SDP program:
		{\small \begin{subequations}\label{sdp_lqr_wconst}
				\begin{align}
					\text{minimize\ \ } &\begin{bmatrix}Q & 0\\0 &R 
					\end{bmatrix} \bullet\Sigma \\
					\text{subject to\ \ } & \Sigma_{xx} = \begin{bmatrix} A_t & B_t\end{bmatrix} \Sigma \begin{bmatrix} A_t & B_t\end{bmatrix}^{\top}+W \label{eq:const_stationary}\\
					& \Sigma \succeq 0 \label{eq:positivedef}\\
					& \begin{bmatrix} A_t & B_t\end{bmatrix} \Sigma \begin{bmatrix} A_t & B_t\end{bmatrix}^{\top} \preceq \alpha \Sigma_{xx} \label{eq:const_stability}
				\end{align}
		\end{subequations}}
		and $K_t = \Sigma_{xu}^{\top} \Sigma_{xx}^{-1}$\;
		\textbf{Play} $u_t = K_t x_t$\;
		\textbf{Update} $x_{t+1} = A_{t} x_t + B_t u_t + w_t, w_t \sim N(0, W)$
	}
	\caption{COCO-LQ: COvariance Constrained Online LQ}
	\label{alg1}
\end{algorithm2e}
COCO-LQ solves an SDP \eqref{sdp_lqr_wconst} at each time step that is similar to that in Proposition~\ref{prop:optimal_lqr_sdp}. 
The crucial difference is the new constraint \eqref{eq:const_stability}, which involves parameter $\alpha$. Plugging  \eqref{eq:const_stationary} into constraint \eqref{eq:const_stability} yields the following:
\[ \Sigma_{xx} \preceq \frac{1}{1-\alpha} W.\]
This highlights that constraint \eqref{eq:const_stability} can be interpreted as an upper bound on the state covariance matrix $\Sigma_{xx}$. When $\alpha = 0$, the controller essentially cancels out the dynamics, without taking into account the cost of doing so.  This ensures stability but can lead to large cost.  At another extreme, when $\alpha \rightarrow 1$, the SDP solved at each time step is the same as for the LTI setting, and so COCO-LQ matches the naive approach in Section~\ref{sec:naive}. Thus, $\alpha$ trades off between stability and cost.
In the following section, we show that this novel state covariance constraint promotes sequential strong stability \citep{cohen2018online}, which in turn guarantees ISS with a proper choice of $\alpha$. 




\vspace{-6pt}
\subsection{Stability}\label{subsec:guarantee}
We now state our main technical result, which provides a formal stability guarantee for COCO-LQ. 
\begin{theorem}\label{thm:stability}
	Let $0 \leq \alpha<1/2$, and suppose \eqref{sdp_lqr_wconst} is feasible for all $t$, then the resulting dynamical system satisfies ISS in the sense that for any disturbance sequence $\{w_t\}_{t=0}^\infty$ and for any $t\geq t_0$, 
	{\small\[\Vert x_t\Vert \leq \rho^{t-t_0} \Vert x_{t_0}\Vert + \frac{\kappa \rho}{1-\rho} \, \sup_{t_0\leq k< t} \Vert w_k\Vert \]}for $\rho=\sqrt{\frac{\alpha}{1-\alpha}}\in [0,1)$ and $\kappa=\frac{\kappa_W}{\sqrt{1-\alpha}}$, where $\kappa_W=\Vert W\Vert \Vert W^{-1}\Vert$ is the condition number of $W$. 
\end{theorem}

\noindent The key intuition underlying this result is that the additional state covariance constraint \eqref{eq:const_stability} implicitly enforces sequential strong stability \citep{cohen2018online}, which in turn ensures ISS. More formally, sequential strong stability is defined as follows,

\begin{definition}{(Sequential Strong Stability)}\label{def:str_seq_stab}
	A sequence of policies $K_1, K_2, ..., $ such that $u_t = K_t x_t$ is $(\kappa, \gamma, \rho)$-sequential strongly stable (for $\kappa > 0$, $0<\gamma\leq 1$ and $0 \leq \rho < 1$) if there exist matrices $H_1, H_2, ..., $ and $L_1,L_2 ..., $ such that $A_t+B_tK_t = H_t L_t H_t^{-1}$ for all $t$, with the following properties: (a)  $\Vert L_t\Vert \leq 1-\gamma$; (b) $\Vert H_t\Vert \leq \beta_1$ and $||H_t^{-1}|| \leq 1/\beta_2$ with $\kappa = \beta_1/{\beta_2}$; (c)  $\Vert H_{t+1}^{-1} H_t\Vert \leq \frac{\rho}{1-\gamma}$.
\end{definition}

\noindent The following lemma formalizes the connection between~\eqref{eq:const_stability} and sequential strong stability.
\setlength{\textfloatsep}{6pt}
\begin{lemma}\label{lemma:seq_stability}
	Under the conditions in Theorem~\ref{thm:stability}, the policies designed by COCO-LQ are $(\kappa, \gamma, \rho)$-sequential strongly stability for $\kappa=\frac{\kappa_W}{\sqrt{1-\alpha}},\gamma=1-\sqrt{\alpha}, \rho=\sqrt{\frac{\alpha}{1-\alpha}}$ where $\kappa_W = \Vert W\Vert \Vert W^{-1}\Vert$. 
\end{lemma}
A proof of Lemma~\ref{lemma:seq_stability} is given in Appendix~\ref{append:proof_maintheorem}. With Lemma~\ref{lemma:seq_stability}, proving the result in Theorem \ref{thm:stability} only requires  showing that sequential strong stability implies ISS. This result is  provided in Appendix~\ref{append:proof_maintheorem} as part of the complete proof of Theorem \ref{thm:stability}.
A critical assumption in Theorem \ref{thm:stability} is that the SDP given in \eqref{sdp_lqr_wconst} is feasible for $0 \leq \alpha <1/2$. The following result shows that when $B_t$ is full row rank, the problem is always feasible. The proof of Lemma~\ref{lemma:feasibility} is postponed to Appendix~\ref{sec:proof_feasibility}. 

\begin{lemma}\label{lemma:feasibility}
	When $B_t$ is full row rank, then the SDP \eqref{sdp_lqr_wconst} is always feasible.
\end{lemma}
\vspace{-6pt}
Note that having $B_t$ full row rank is a sufficient but not necessary condition for feasibility of \eqref{sdp_lqr_wconst} of COCO-LQ. When $B_t$ is not full row rank, the feasibility assumption may still hold, and therefore our assumption is weaker than the invertibility assumption used in the literature, e.g.  \citet{LAI198623}. 
\revise{More broadly, in Theorem~\ref{thm:stability}, $\alpha<0.5$ is a sufficient condition for stability. For $\alpha > 0.5$, stability may still hold for some problem instances $(A_t, B_t)$ as will be shown in the simulations in Section~\ref{sec:experiments}. How to provide a more refined instance-dependent threshold on $\alpha$ is an interesting future direction.}

\vspace{-12pt}
\subsection{Infeasibility and the Role of Predictions}\label{subsec:infeasibility}
We now turn our attention to the case when the SDP given in \eqref{sdp_lqr_wconst} is infeasible.  In this case it is necessary for the controller to use additional information in order to stabilize the system. In particular, we provide an example that shows the necessity of predictions when 
$B_t$ is not full row rank in Appendix~\ref{append:example_pred}. This example shows that when $B_t$ is not full row rank, for any (deterministic) online control algorithm that has causal access to system matrices, there exists a future sequence of $(A_t, B_t)$ in which the algorithm cannot stabilize the system. 
\setlength{\textfloatsep}{6pt}
\begin{algorithm2e}[t]
	\SetAlgoLined
	\DontPrintSemicolon
	\LinesNumbered
	\Parameter{$\alpha \in [0, 1)$}
	\KwIn{$Q, R, W \succ 0$}
	\For{$t=1, 2, ...$}{
		\If{ $t \equiv 1 \ (\mathrm{mod}\ H)$}{
			\textbf{Receive} state $x_t$, and system parameters $(A_t, B_t), \ldots, (A_{t+H-1}, B_{t+H-1})$\;
			\textbf{Compute policy:} Let $\Sigma_{t} \in \mathbb{R}^{n \times n}$ be a solution to the constrained SDP in ~\eqref{sdp_lqr_wconst} with
			$(R, A_t, B_t)$ replaced by $(\tilde{R}, \tilde{A}_t, \tilde{B}_t)$, where $\tilde{R} = \begin{bmatrix}
			R & & \\
			& \ddots & \\
			& & R
			\end{bmatrix}$ ($H$ repeating blocks), 
			$\tilde{A}_t:= A_{t+H-1}\cdots A_t, \tilde{B}_t:= [B_{t+H-1}, A_{t+H-1}B_{t+H-2}, \cdots,  A_{t+H-1}...A_{t+1} B_t]$ \;
			Set $K_t = \Sigma_{xu}^{\top} \Sigma_{xx}^{-1}$ and $[u_{t+H-1}, ..., u_t]^\top = K_t x_t$\;
		}
		\textbf{Play} Implement the planned control action $u_t$\;
		\textbf{Update} $x_{t+1} = A_{t} x_t + B_t u_t + w_t, w_t \sim \mathcal{N}(0, W)$
	}
	\caption{COCO-LQ-Prediction: COVariance Constrained Online LQ with Predictions}
	\label{alg2}
\end{algorithm2e}
In this section, we show that using $(A_t, B_t)$ together with short-term predictions of future system matrices is enough to stabilize the system under standard controllability assumptions. Specifically, we extend COCO-LQ to include future $H$ steps of predictions in Algorithm \ref{alg2}. The key idea is to rewrite the dynamics as
\begin{align}
	x_{t+H} = \tilde{A}_t x_t + \tilde{B}_t \bar{u}_t 
	+ [I, A_{t+H-1}, \cdots,  A_{t+H-1}...A_{t+1}] \bar{w}_t 
\end{align}
where we define $\tilde{A_t}:= A_{t+H-1}\cdots A_t$, $\tilde{B_t} := [B_{t+H-1}, A_{t+H-1}B_{t+H-2}, \cdots,  A_{t+H-1}...A_{t+1} B_t]$, $\bar{u}_t \!:=\! [u_{t+H-1}^\top, u_{t+H-2}^\top, \ldots, u_t^\top]^\top$ and $\bar{w}_t \!:=\! [w_{t+H-1}^\top, w_{t+H-2}^\top, \ldots, w_t^\top]^\top$. When $H$ is long enough such that $\tilde{B}_t$ is full row rank, we can use Algorithm~\ref{alg1} on $\tilde{A_t}$ and $\tilde{B_t}$ and avoid the infeasibility issue, and our stability guarantee is provided below. The proof of Theorem~\ref{thm:predictions} is in Appendix~\ref{sec:predictions}.

\begin{theorem}\label{thm:predictions}
	Suppose for each $t$, matrix $\tilde{B}_t = [B_{t+H-1}, A_{t+H-1}B_{t+H-2}, \cdots,  A_{t+H-1}...A_{t+1} B_t]$ satisfies $\tilde{B}_t \tilde{B}_t^\top \succeq \sigma I $  for some $\sigma>0$, and $\|A_t\|\leq a, \|B_t\|\leq b$ for some $a,b>0$. Then, the SDP in Algorithm~\ref{alg2} is always feasible. Further, when $\alpha<1/2$, the closed-loop system is ISS for any $t$,
	{\small\begin{align*}
			\Vert x_{t}\Vert  \leq \kappa_A'  \rho^{\frac{t}{H} -1 }\Vert x_{1}\Vert +\kappa_A'\kappa_A \kappa \max(1, \frac{\rho}{1-\rho}) \sup_{1 \leq s <  t} \Vert w_s\Vert,
	\end{align*}}
	where the same as Theorem~\ref{thm:stability}, $\rho=\sqrt{\frac{\alpha}{1-\alpha}}\in [0,1)$ and $\kappa=\frac{\kappa_W}{\sqrt{1-\alpha}}$ with $\kappa_W=\Vert W\Vert \Vert W^{-1}\Vert$ being the condition number of $W$; further, $\kappa_A = 1+a+\ldots+a^{H-1}$, and $\kappa_A' = a^{H-1} + b^2(1+a+\cdots + a^{H-1})^2 \kappa_R \frac{ \kappa + a^{H}}{\sigma}$ with $\kappa_R$ being the condition number of $R$. 
\end{theorem}
\vspace{-12pt}
\revise{\subsection{Estimation Error}
	\label{sec:estimation_error}
	In both Algorithm~\ref{alg1} and Algorithm~\ref{alg2}, the exact knowledge of state-transition matrices $(A_t, B_t)$ or the extended state-transition matrices $(\tilde{A}_t, \tilde{B}_t)$ are needed when deriving the control actions. 
	In this section, we show that COCO-LQ can still obtain a stabilizing controller in the case where only approximations are known, if the estimation error is controlled. Our main result is the following.
	
	\begin{theorem}\label{thm:estimation}
		Let $(\hat{A}_t$, $\hat{B}_t)$ be an estimate of ($A_t$, $B_t$). Given $\alpha \in [0, \frac{1}{2})$, let $\rho = \sqrt{\frac{\alpha}{1-\alpha}}$, $\kappa = \frac{\Vert W\Vert \Vert W^{-1}\Vert}{\sqrt{1-\alpha}}$ and $\gamma = 1-\sqrt{\alpha}$. Let $K_1, K_2, ...$ be the policies designed by COCO-LQ for $(\hat{A}_t$, $\hat{B}_t)$ with parameter $\alpha$. When the estimation error satisfies,
		{\small\begin{equation}
			\max\{||\hat{A}_t-A_t||_2, ||\hat{B}_t -B_t||_2\} \leq \delta \frac{\gamma}{\kappa(1+K_{max})} \label{eq:estimation_err}
			\end{equation}}where $\delta $ can be any number in $(0,\frac{\sqrt{1-\alpha}-\sqrt{\alpha}}{1-\sqrt{\alpha}})$, and $K_{max}$ is any uniform upper bound on $\Vert K_t\Vert$.
		Then, the policies $K_t$ are ISS when applied to the system $(A_t, B_t)$, 
		{\small\[\Vert x_t\Vert \leq \left(\rho'\right)^{t-t_0} \Vert x_{t_0}\Vert + \frac{\kappa\rho'}{1-\rho'}\sup_{t_0\leq k<t} \Vert w_k\Vert, \]
			where $\rho' = \frac{1-(1-\delta) \gamma}{1-\gamma} \rho \in (0,1)$. Finally, when $\Vert \hat{A}_t\Vert  \leq \bar{\sigma}_A$, $  \Vert \hat{B}_t\Vert  \leq \bar{\sigma}_B$ and $\hat{B}_t \hat{B}_t^\top \succeq \underline{\sigma}_B^2$, one uniform upper bound for $\Vert K_t\Vert$ is $K_{max} = \kappa_R \frac{\bar{\sigma}_B}{\underline{\sigma}_B^2}(\kappa (1-\gamma) + \bar{\sigma}_A)$ with $\kappa_R = \Vert R\Vert \Vert R^{-1}\Vert$. }
	\end{theorem}
	A proof of Theorem~\ref{thm:estimation} is provided in Appendix~\ref{sec:estimation}. This result highlights the tradeoff between the estimation error and the algorithm performance. If we choose a small $\alpha$, the algorithm can tolerant a larger estimation error (i.e. larger right hand side of \eqref{eq:estimation_err} can be obtained) but may lead to high control cost due to the tight state co-variance constraint. If we choose a larger $\alpha$, the algorithm tolerates smaller estimation error while its performance improves due to the less strict state co-variance constraint. 
	
}

\section{Experiments}\label{sec:experiments}
The results in the previous section focus on stability of COCO-LQ approach. Here, we use experimental results to highlight that COCO-LQ also performs near-optimally in terms of cost while also stabilizing systems that the naive approach based on LTI control cannot. 
In Section~\ref{sec:synthetic}, we test our method on random, synthetic linear time varying systems, and in Section~\ref{sec:power} we demonstrate the algorithm performance in real-world power system frequency control settings.
More experiment results on nonlinear systems via local linear approximation could be found in Appendix~\ref{sec:nonlinear}. 

\subsection{Synthetic Time-Varying Systems}
\label{sec:synthetic}
We first consider the control of switching and time-variant systems. The cost function is set as $Q = 0.2I, R = I$, and system is subject to Gaussian disturbance $w_t \sim N(0, 0.1^2)$. We average the simulation results over 5 runs and visualize the mean performance and standard deviation.
\begin{enumerate}[leftmargin=*, label=\alph*.]
\item \textit{Switching systems.} we consider a switching system following Example~\ref{example} in Section~\ref{sec:naive}, where $A_t$ alternates between $A = [[0.99, 1.5], [0, 0.99]]$ and $A' = [[0.99, 0], [1.5, 0.99]]$, and $B_t = I$.
\item \textit{Time-variant systems.} We consider a system $A_t = [[0.99, \sin(\frac{\pi t}{2})|e^{t/60}], [|\cos(\frac{\pi t}{2})|e^{t/60}, 0.99]]$ that is continually changing over time, and $B_t = I$.
\end{enumerate}
\begin{figure}[htbp]
\includegraphics[width=0.24\columnwidth]{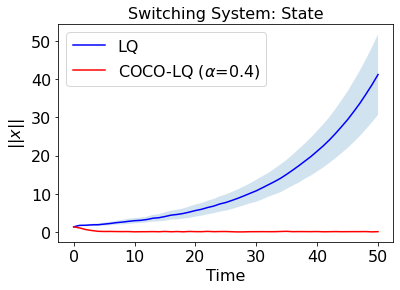}
\includegraphics[width=0.24\columnwidth]{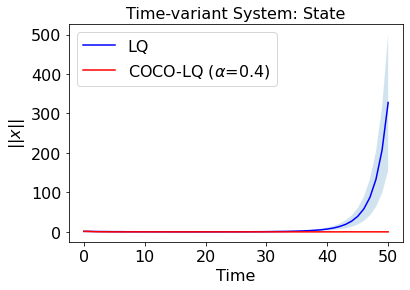}
\includegraphics[width=0.24\columnwidth]{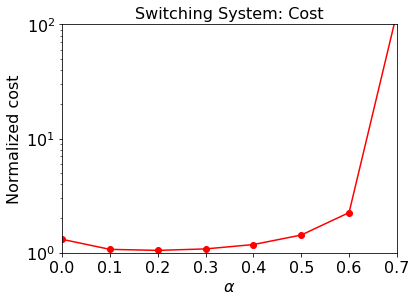}
\includegraphics[width=0.24\columnwidth]{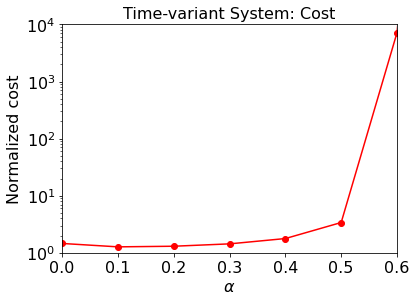}
\caption{Performance comparison of COCO-LQ and LQ on synthetic time-varying systems. The left two figures show the state evolution, and right two figures show the normalized cost (cost of COCO-LQ divided by cost of the offline optima) under different $\alpha$.}
\label{fig:synthetic_result}
\end{figure}
As we can see in Figure~\ref{fig:synthetic_result}, COCO-LQ is able to quickly and effectively stabilize the system under various time-varying scenarios, which validates our theoretical findings. As $\alpha$ increases, the acquired cost of COCO-LQ first decreases and then increases (explosion of state), highlighting that $\alpha$ can explicitly control the tradeoff between cost and stability. With proper selection of $\alpha$, COCO-LQ achieves near-optimal cost (within 30\% of the offline optimal for both system a and b).



\subsection{Frequency Control with Renewable Generation}
\label{sec:power}
\begin{figure}[t]
\begin{subfigure}
  \centering
  \includegraphics[width=0.3\columnwidth]{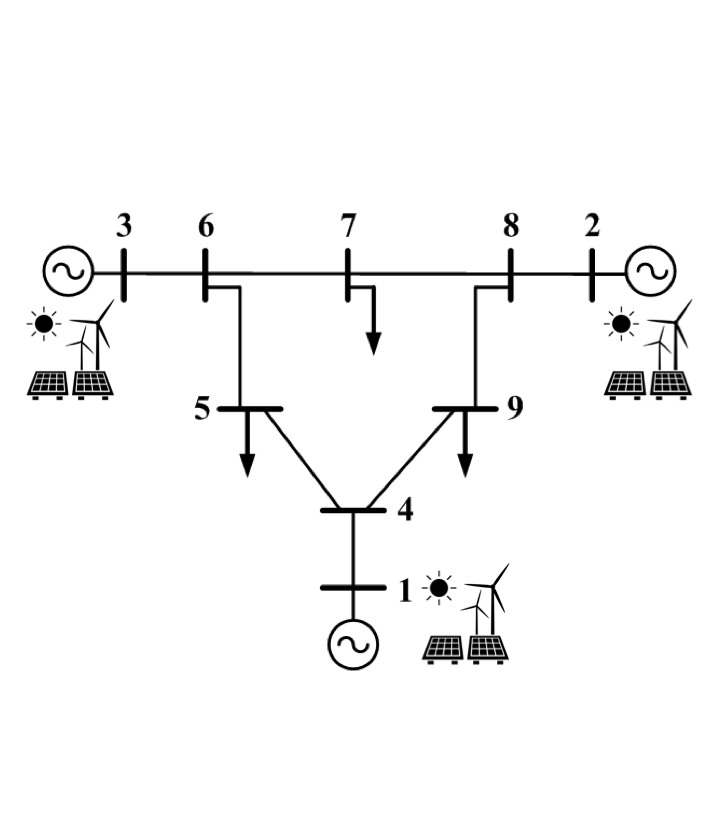}
  \label{fig:ieee_case}
\end{subfigure}%
\begin{subfigure}
  \centering
  \includegraphics[width=0.6\columnwidth]{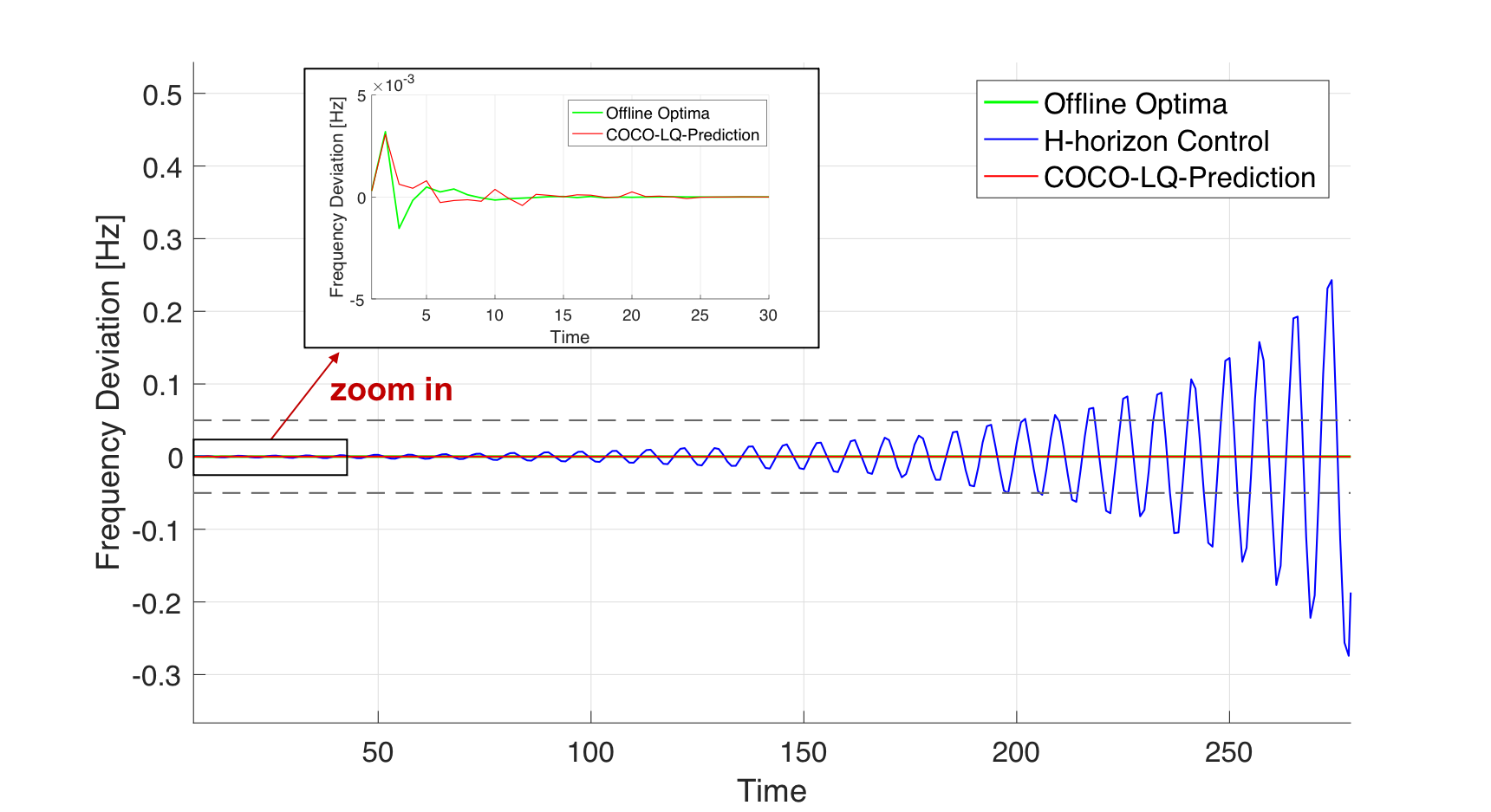}
  \label{fig:results}
\end{subfigure}
\caption{(a) IEEE WECC 3-machine 9-bus system schematic with generators at bus 1, 5, 9 are mixture of thermal generation and renewable. (b) Frequency dynamics under offline optima, baseline H-horizon control, and COCO-LQ. The dotted grey lines ($\pm 0.05$Hz) are the safety margin of power system frequency variation.}
\label{fig:frequency_dyn}
\end{figure}

We now consider a power system frequency control problem on standard IEEE WECC 3-machine
9-bus system (Figure~\ref{fig:frequency_dyn}(a)), which is a widely adopted system used in frequency stability studies. The state space model of power system frequency dynamics follows~\cite{hidalgo2019frequency},
\begin{equation}\label{eq:freq_dyn}
    \underbrace{\begin{bmatrix}
    \dot{\theta}\\\dot{\omega} 
    \end{bmatrix}}_{\dot{\mathbf{x}}}=
    \underbrace{\begin{bmatrix}
    0 & I \\
    -M_t^{-1}L & -M_t^{-1}D
    \end{bmatrix}}_{A_t}
  \begin{bmatrix}
    \theta\\\omega
    \end{bmatrix}
    + \underbrace{\begin{bmatrix}
    0\\M_t^{-1}
    \end{bmatrix}}_{B_t} \underbrace{p_{in}}_{\mathbf{u_t}}
\end{equation}
where the state variable is defined as the stacked vector of the voltage angle $\theta$ and frequency $\omega$. $M_t = diag(m_{t, i})$ is the inertia matrix, where $m_{t, i}$ represents the equivalent rotational inertia at bus $i$ and time $t$. $M_t$ is time-varying and depends on the mix of online generators, since only thermal generators provide rotational inertia and renewable generation does not~\cite{ulbig2014impact}. $D= diag(d_i)$ is the damping matrix, where $d_i$ is the generator damping coefficient. $L$ is the network susceptance matrix. The control variable $p_{in}$ corresponds to the electric power generation.

We assume the system is changing between two states: a high renewable generation scenario where $m_{t, i}=2$ (i.e., 80 percent renewable with zero inertia and 20 percent of thermal generation with 10s inertia), and a low renewable generation scenario where $m_{t, i}=8$ (i.e., 20 percent renewable and 80 percent thermal generation), with additional random fluctuations between $[0, 0.2]$. This setup represents the real-world situation where we have high solar output
during the daytime, and low output in the morning/evening, with intra-day variations due to clouds and weather changes. Notice that $B_t$ is not full rank, thus we need to leverage predictions, i.e., $A_{t+1}$ and $B_{t+1}$. For fair comparison, we compete against the $H$-horizon optimal control in~\cite{bertsekas1995dynamic}, which is the extension of naive LTI controller to use $H$-step predictions. In both cases, we assume the prediction is accurate and use the exact value of $A_{t+1}$ and $B_{t+1}$ for computing control actions. 

Figure~\ref{fig:frequency_dyn}(b) visualizes the power system frequency dynamics under three controllers: the offline optimal control, the baseline $H$-horizon optimal controller and the proposed COCO-LQ-Prediction method. We ideally desire a controller that is able to maintain the frequency variation within $\pm 0.05$Hz and eventually stabilize the system. It can be observed that our algorithm succeeds at maintaining the frequency stability under random, time-varying renewable generations. Furthermore, the performance of COCO-LQ is very close to the offline optimal, while the system frequency diverges under the baseline $H$-horizon optimal control. 

\section{Conclusion}
\revise{In this paper, we study the stability of LTV systems. Our results demonstrate the challenge of ensuring stability for LTV systems compared to LTI systems.  Motivated by this challenge, we propose a COCO-LQ/COCO-LQ-Prediction policy that can guarantee stability for LTV systems under certain assumptions. There are many interesting open questions that remain. For example, the bound $\alpha<1/2$ in Theorem~\ref{thm:stability} is a sufficient condition, and studying how to relax the bound and how to derive instance-dependent bounds is an interesting future question. Another important direction is to analyze the performance (e.g. the regret) of the proposed approach in order to quantify the tradeoff between stability and performance.  
}






\bibliography{reference.bib}

\newpage
\appendix
\section{Proof of Theorem~\ref{thm:stability}}
\label{append:proof_maintheorem}
In this section, we provide the proof of Theorem~\ref{thm:stability}. As mentioned in Section~\ref{subsec:guarantee}, the proof consists of two steps. In the first step, we show that the algorithm guarantees strong sequential stability (Lemma~\ref{lemma:seq_stability}), which we prove in Appendix~\ref{subsec:str_seq}, and in the second step, we show that strong sequential stability guarantees input-to-state stability, or more formally:
\begin{lemma}\label{lem:iss}
Supppose a sequence of policies $K_1,K_2,\ldots,K_t,\ldots$ is $(\kappa,\gamma,\rho)$-sequential strongly stable. Then, the closed loop system is input-to-state stable in the sense that for any $t\geq t_0\geq 1$,
\begin{align*}
    \| x_t \|  \leq \kappa \rho^{t-t_0} \|x_{t_0}\| + \frac{\kappa \rho}{1-\rho} \max_{t_0 \leq s <t}\|w_s\|.
\end{align*}
\end{lemma}
We provide the proof of Lemma~\ref{lem:iss} in Appendix~\ref{subsec:iss}. Combining the above two steps finishes the proof of Theorem~\ref{thm:stability}. 

\subsection{Proof of Lemma~\ref{lemma:seq_stability}}\label{subsec:str_seq}
To prove Lemma~\ref{lemma:seq_stability}, we first show an auxilliary result that the optimal solution to the SDP in \eqref{sdp_lqr_wconst} must have a specific low rank structure. Then, utilizing this structure, we prove Lemma~\ref{lemma:seq_stability}. 
\begin{lemma}\label{lem:realizability}
Suppose $Q,R,W$ are strictly positive definite and $\alpha\in[0,1)$. If \eqref{sdp_lqr_wconst} has a minimizer $\Sigma^*$, then there exists $K\in\mathbb{R}^{p\times d}$ s.t. $\Sigma^*$ can be written as  
\begin{align*}
 \Sigma^*= \begin{bmatrix}\Sigma_{xx}^{*} & \Sigma_{xu}^{*}\\
 {(\Sigma_{xu}^*)^{\top}} & \Sigma_{uu}^{*}\end{bmatrix} = \begin{bmatrix}\Sigma_{xx}^{*} & \Sigma_{xx}^{*} K^{\top}\\
K \Sigma_{xx}^{*} & K \Sigma_{xx}^{*} K^{\top} \end{bmatrix}
\end{align*}
\end{lemma}
\begin{proof}
Clearly, $\Sigma_{xx}^* \succeq W \succ 0$, and therefore we can simply define $K = (\Sigma_{xu}^*)^\top(\Sigma_{xx}^*)^{-1}$, and the only thing we need to show is $\Sigma_{uu}^* = K \Sigma_{xx}^{*} K^{\top}$. Suppose this is not true, then since $\Sigma^*\succeq 0$, there must exist $D\neq 0, D\succeq 0$ s.t. 
\begin{align}
    \Sigma_{uu}^* = K \Sigma_{xx}^{*} K^{\top} + D.
\end{align}
Then, by \eqref{eq:const_stationary}, we also have,
\begin{align}
    \Sigma_{xx}^* = (A_t + B_t K)\Sigma_{xx}^* (A_t+ B_tK)^{\top} + B_t D B_t^{\top}  + W. \label{eq:realizability:sigmaxx}
\end{align}  
Viewing the above as a Lyapunov equation in terms of $\Sigma_{xx}^*$, and since $B_t D B_t^{\top} + W\succ 0$ and $\Sigma_{xx}^*\succ 0$, we get $A_t + B_t K$ is a stable matrix. 

Next, since $A_t + B_t K$ is stable, we can construct $\tilde{\Sigma}_{xx}$ to be the unique solution to the following Lyapunov equation,
\begin{align}
    \tilde\Sigma_{xx} = (A_t + B_t K)\tilde\Sigma_{xx} (A_t+ B_tK)^{\top}  + W. \label{eq:realizability:tildesigmaxx} 
\end{align}
We further define,
\begin{align*}
 \tilde{\Sigma}=\begin{bmatrix}\tilde\Sigma_{xx} & \tilde\Sigma_{xx} K^{\top}\\
K \tilde\Sigma_{xx} & K \tilde\Sigma_{xx} K^{\top} \end{bmatrix},
\end{align*}
and we claim that $\tilde\Sigma$ is a feasible solution to \eqref{sdp_lqr_wconst}. Clearly $\tilde{\Sigma}$ is positive semi-definite, and further, it satisfies \eqref{eq:const_stationary}. We now check that \eqref{eq:const_stability} is met below. 
\begin{align}
    &\begin{bmatrix} A_t & B_t\end{bmatrix} \tilde\Sigma \begin{bmatrix} A_t & B_t\end{bmatrix}^\top - \alpha \tilde{\Sigma}_{xx} \nonumber \\
    &=( A_t + B_t K)\tilde{\Sigma}_{xx} (A_t + B_tK)^\top - \alpha \tilde{\Sigma}_{xx}\nonumber\\
    &= ( A_t + B_t K){\Sigma}_{xx}^* (A_t + B_tK)^\top - \alpha {\Sigma}_{xx}^* + ( A_t + B_t K) ( \tilde{\Sigma}_{xx} -  {\Sigma}_{xx}^* ) (A_t + B_tK)^\top - \alpha( \tilde{\Sigma}_{xx} -  {\Sigma}_{xx}^*)\nonumber\\
    &= \begin{bmatrix} A_t & B_t\end{bmatrix} \Sigma^* \begin{bmatrix} A_t & B_t\end{bmatrix}^\top - B_t D B_t^\top - \alpha {\Sigma}_{xx}^* + ( A_t + B_t K) ( \tilde{\Sigma}_{xx} -  {\Sigma}_{xx}^* ) (A_t + B_tK)^\top - \alpha( \tilde{\Sigma}_{xx} -  {\Sigma}_{xx}^*)\nonumber \\
    &\preceq - B_t D B_t^\top + ( A_t + B_t K) ( \tilde{\Sigma}_{xx} -  {\Sigma}_{xx}^* ) (A_t + B_tK)^\top - \alpha( \tilde{\Sigma}_{xx} -  {\Sigma}_{xx}^*), \label{eq:realizability:check4d}\end{align}
    where the last step is due to $\Sigma^*$ must satisfy \eqref{eq:const_stability}. 
    Subtracting \eqref{eq:realizability:sigmaxx} from \eqref{eq:realizability:tildesigmaxx}, we have,
    \begin{align}
       \tilde{\Sigma}_{xx} - \Sigma_{xx}^* =  ( A_t + B_t K) ( \tilde{\Sigma}_{xx} -  {\Sigma}_{xx}^* ) (A_t + B_tK)^\top  - B_t D B_t^\top.  \label{eq:realizability:check4d_2}
    \end{align}
    Plugging \eqref{eq:realizability:check4d_2} into \eqref{eq:realizability:check4d}, we have, 
  \begin{align*}
          \begin{bmatrix} A_t & B_t\end{bmatrix} \tilde\Sigma \begin{bmatrix} A_t & B_t\end{bmatrix}^\top - \alpha \tilde{\Sigma}_{xx} 
    \preceq (1-\alpha) (\tilde{\Sigma}_{xx} - \Sigma_{xx}^*) .
  \end{align*}
 Therefore, to check \eqref{eq:const_stability} it remains to show $\tilde{\Sigma}_{xx} - \Sigma_{xx}^*\preceq 0$. To see this, we view \eqref{eq:realizability:check4d_2} as a Lyapunov equation in terms of $\tilde{\Sigma}_{xx} - \Sigma_{xx}^*$, and as $A_t + B_tK$ is stable, and $B_t D B_t^\top \succeq  0$, we have $\tilde{\Sigma}_{xx} - \Sigma_{xx}^* \preceq 0$. As a result, \eqref{eq:const_stability} holds and 
  $\tilde{\Sigma}_{xx} $ is indeed a feasible solution to \eqref{sdp_lqr_wconst}. 
  
  Further, we can show $\tilde{\Sigma} $ achieves a strictly lower cost than $\Sigma^*$. 
  To see this, note we have already shown $\tilde{\Sigma}_{xx} \preceq \Sigma_{xx}^*$, which also implies  $\tilde{\Sigma}_{uu} = K \tilde{\Sigma}_{xx} K^\top \preceq K \Sigma^*_{xx} K^\top =  \Sigma_{uu}^* - D$ with $D\succeq 0$, $D\neq 0$. Coupled this with the fact that $Q,R$ are strictly positive definite, we have $\tilde{\Sigma}$ must achieve a lower cost. So we get a contradiction, and the proof is concluded. 
\end{proof}

Leveraging the decomposition of the SDP solution in Lemma~\ref{lem:realizability}, we now prove Lemma~\ref{lemma:seq_stability}. 

\noindent\textbf{Proof of Lemma~\ref{lemma:seq_stability}}. Let $\Sigma_t$ be the solution to the SDP \eqref{sdp_lqr_wconst} at step $t$. Then by Lemma~\ref{lem:realizability}, $\Sigma_t$ can be rewritten as,
\[\Sigma_t = \begin{bmatrix}\Sigma_{t,xx} & \Sigma_{t,xx}K_t^{\top}\\ K_t\Sigma_{t,xx} & K_t\Sigma_{t,xx}K_t^{\top}\end{bmatrix},  \]
for some $\Sigma_{t,xx}\succ 0$, and $K_t$ is the linear controller at step $t$. With this, we can re-write the left side of constraint~\eqref{eq:const_stability} as follows,
\begin{align*}
    \begin{bmatrix} A_t & B_t\end{bmatrix} \Sigma_t \begin{bmatrix} A_t & B_t\end{bmatrix}^{\top} &= \begin{bmatrix} A_t & B_t\end{bmatrix} \begin{bmatrix}\Sigma_{t,xx} & \Sigma_{t,xx}K_t^{\top}\\ K_t\Sigma_{t,xx} & K_t\Sigma_{t,xx}K_t^{\top}\end{bmatrix} \begin{bmatrix} A_t^{\top} \\ B_t^{\top}\end{bmatrix} \\
    &= A_t \Sigma_{t,xx} A_t^{\top} + A_t K_t\Sigma_{t,xx} B_t^{\top}+B_t\Sigma_{t,xx} K_t^{\top} A_t^{\top} + B_t K_t \Sigma_{t,xx} K_t^{\top} B_t^{\top} \\
    &= (A_t+B_tK_t) \Sigma_{t,xx} (A_t+B_t K_t)^{\top}.
\end{align*}
As a result, \eqref{eq:const_stability} can be equivalently expressed as,
\begin{equation}
    (A_t+B_tK_t) \Sigma_{t,xx} (A_t+B_t K_t)^{\top} \preceq \alpha \Sigma_{t,xx}.
\end{equation}
Left and right multiplying the above by $\Sigma_{t,xx}^{-1/2}$, we get
\begin{equation}\label{eq:lt_bound}
\Sigma_{t,xx}^{-1/2} (A_t+B_tK_t) \Sigma_{t,xx}^{1/2} \Sigma_{t,xx}^{1/2} (A_t+B_t K_t)^{\top} \Sigma_{t,xx}^{-1/2} \preceq \alpha I\,.
\end{equation}
Let $H_t = \Sigma_{t,xx}^{1/2}$, $L_t = \Sigma_{t,xx}^{-1/2} (A_t+B_tK_t) \Sigma_{t,xx}^{1/2}$. This gives us the following decomposition,
\begin{equation}
    (A_t+B_tK_t) = H_t L_t H_t^{-1}\,.
\end{equation}
We now show that this decomposition yields the desired sequential stability property. To do so, we need to check the three conditions in Definition~\ref{def:str_seq_stab}. 

To check condition (a) in Definition~\ref{def:str_seq_stab}, note that
\eqref{eq:lt_bound} provides an upper bound on $\|L_t\|$, that is $L_t L_t^{\top} \preceq \alpha I$. Therefore, $\|L_t\| \leq \sqrt{\alpha} = 1-\gamma$ with $\gamma = 1-\sqrt{\alpha}$.  

To check condition (b) which is an upper bound on $\|H_t\|$ and $ \|H_t^{-1}\|$, we recall constraint \eqref{eq:const_stationary}, 
$$\begin{bmatrix} A_t & B_t\end{bmatrix} \Sigma_t \begin{bmatrix} A_t & B_t\end{bmatrix}^{\top} = \Sigma_{t,xx}-W.$$
As the left hand side of the above is positive semi-definite, we must have $\Sigma_{t,xx}\succeq W$. Further, plugging
\eqref{eq:const_stationary} into \eqref{eq:const_stability}, we can get,
$\Sigma_{t,xx}- W \preceq \alpha \Sigma_{t,xx}$. These can be equivalently written as, 
\begin{align}
    W \preceq \Sigma_{t,xx} \preceq \frac{W}{1-\alpha}. \label{eq:str_seq:sigma_ublb}
\end{align}

Therefore, $\|H_t\| = \|\Sigma_{t,xx}^{1/2}\| \leq \frac{\|W^{1/2}\|}{\sqrt{1-\alpha}}$, $\|H_t^{-1}\| = \|\Sigma_{t,xx}^{-1/2}\| \leq {\|W^{-1/2}\|}$, $\|H_t\| \|H_t^{-1}\| \leq \frac{\| W^{1/2}\| \| W^{-1/2}\|}{\sqrt{1-\alpha}} = \frac{\kappa_W}{\sqrt{1-\alpha}}$. Therefore, condition (b) holds with $\kappa = \frac{\kappa_W}{\sqrt{1-\alpha}}$. 

Lastly, we check condition (c), which is an upper bound on $\Vert H_{t+1}^{-1} H_t\Vert$. Note that,
\begin{align*}
    H_{t+1}^{-1} H_t H_t^\top (H_{t+1}^{-1})^\top= H_{t+1}^{-1}  \Sigma_{t,xx} (H_{t+1}^{-1})^\top \preceq   \frac{1}{1-\alpha}(\Sigma_{t+1,xx})^{-1/2}  W (\Sigma_{t+1,xx})^{-1/2}\preceq \frac{1}{1-\alpha} I,
\end{align*}
where the two inequalities in the above derivations are due to \eqref{eq:str_seq:sigma_ublb}. As a result, we have $\|H_{t+1}^{-1} H_t\| \leq \frac{1}{\sqrt{1-\alpha}} = \frac{\rho}{1-\gamma}$ with $\rho = \frac{\sqrt{\alpha}}{\sqrt{1-\alpha}}$.  As $\alpha<1/2$, we have $\rho<1$. As such, condition (c) holds.

In summary, the $(\kappa,\gamma,\rho)$ strong sequential stability holds, with $\kappa = \frac{\kappa_W}{\sqrt{1-\alpha}}$, $\gamma = 1-\sqrt{\alpha}$, and $\rho = \sqrt{\frac{\alpha}{1-\alpha}}$. This concludes the proof.
\qed

\subsection{Proof of Lemma~\ref{lem:iss}} \label{subsec:iss}
For notational simplicity, we only prove the case with $t_0 = 1$; the general case follows similarly, at the cost of heavier notations. Let $x_1, x_2, ...$ be a sequence of states starting from initial state $x_1$, and generated by the dynamics $x_{t+1} = A_t x_t + B_t u_t + w_t= (A_t+B_t K_t)x_t + w_t$. Hence,
$$x_t = M_1 x_1 + \sum_{s=1}^{t-1} M_{s+1} w_s\,,$$
where 
$$M_t = I; M_s=M_{s+1}(A_s+B_s K_s) = (A_{t-1}+B_{t-1}K_{t-1}) \cdots (A_{s}+B_s K_s). 
$$

By the strong sequential stability, there exist matrices $H_1, H_2,...$ and $L_1, L_2, ...$ such that $A_j + B_j K_j = H_j L_j H_j^{-1}$, and $H_j, L_j$ satisfy the properties specified in Definition~\ref{def:str_seq_stab}. Thus we have for all $1 \leq s < t$,
\begin{align*}
    \| M_s \| &= \| H_{t-1} L_{t-1} H_{t-1}^{-1} H_{t-2} L_{t-1} H_{t-2}^{-1} \cdots H_s L_s H_s^{-1} \| \\
    & \leq \|H_{t-1}\| \left(\prod_{j=s}^{t-1} \|L_{j}\|\right) \left(\prod_{j=s}^{t-2} \|H_{j+1}^{-1} H_j\| \right) \|H_s^{-1}\|\,,\\
    & \leq \beta_1 (1-\gamma)^{t-s} (\frac{\rho}{1-\gamma})^{t-s-1} (1/\beta_2) \leq \frac{\kappa(1-\gamma)}{\rho} \rho^{t-s} 
\end{align*}
As $\kappa \geq 1$, the same holds for $M_t$. Thus, we have, 
\begin{align*}
    \| x_t \| &\leq \| M_1 \|  \| x_1\| + \sum_{s=1}^{t-1} \| M_{s+1}\| \|w_s\| \\
    & \leq \kappa \rho^{t-1} \|x_1\| + \kappa \sum_{s=1}^{t-1} \rho^{t-s-1} \|w_s\| \\
    & \leq \kappa \rho^{t-1} \|x_1\| + \kappa \max_{1 \leq s <t}\|w_s\| \sum_{t=1}^{\infty} \rho^{t}\\
    & = \kappa \rho^{t-1} \|x_1\| + \frac{\kappa \rho}{1-\rho} \max_{1 \leq s <t}\|w_s\|.
\end{align*}
 
 \section{Proof of Lemma~\ref{lemma:feasibility}}\label{sec:proof_feasibility}
 
We prove the lemma by explicitly constructing a feasible solution for ~\eqref{sdp_lqr_wconst}. As $B_t$ has full row rank, $B_t B_t^\top$ is invertible. 
Let \begin{align*}\Sigma_0 &= \begin{bmatrix}W & -WA_t^{\top} (B_t B_t^\top)^{-1} B_t  \\
-B_t^\top (B_t B_t^\top)^{-1}A_t W & B_t^\top (B_t B_t^\top)^{-1}A_t W A_t^{\top} (B_t B_t^\top)^{-1} B_t \end{bmatrix} \\
&= \begin{bmatrix}
I \\
-B_t^\top (B_t B_t^\top)^{-1}A_t
\end{bmatrix} W \begin{bmatrix}
I \\
-B_t^\top (B_t B_t^\top)^{-1}A_t
\end{bmatrix}^\top. \end{align*} 
It suffices to show that $\Sigma_0$ satisfies~\eqref{eq:const_stationary},~\eqref{eq:positivedef} and~\eqref{eq:const_stability}. Notice that
\begin{align*}
    \begin{bmatrix} A_t & B_t\end{bmatrix} \Sigma_0 \begin{bmatrix} A_t & B_t\end{bmatrix}^{\top} = (A_t - B_t B_t^\top (B_t B_t^\top)^{-1}A_t )W (A_t - B_t B_t^\top (B_t B_t^\top)^{-1}A_t )^\top = 0.
\end{align*}
Therefore, constraint \eqref{eq:const_stationary} is equivalent to $\Sigma_{xx} = W$, which holds by the construction of $\Sigma_0$. Next, as $W\succ 0$, $\Sigma_0$ is clearly positive semi-definitive and \eqref{eq:positivedef} holds. Finally, note the left hand side of \eqref{eq:const_stability} is $0$, and the right hand side of \eqref{eq:const_stability} is positive semi-definite. As a result, \eqref{eq:const_stability} holds.  

\section{On the Role of Prediction}
\label{append:example_pred}
In the following example, we show that when $B_t$ is not full row rank, then no matter what (deterministic) online algorithm we use, there always exists a time varying $(A_t,B_t)$ for which this algorithm fails to stabilize from a initial state. 

\begin{example} We consider a counter example which shows when $B_t$ is not full row rank, then no matter what online deterministic algorithm we use which decides $u_t$ based on information of the dynamics up to $(A_t,B_t)$, we can construct a sequence of $(A_t,B_t)$ such that the system state will blow up from a given initial state. 

We set $d=2,p=1$, i.e. $A_t$ is $2$-by-$2$ and $B_t$ is $2$-by-$1$. For all $t$, we set $B_t = [1,0]^\top$, and for even time indices $t = 2k$, we set $A_{t} = I$. Further, $x_0 = [1,1]^\top$ and $w_t = 0$. Our construction is based on induction. Suppose at an even index $t=2k$, $x_{2k} = [x_{2k,1},x_{2k,2}]^\top$ with $x_{2k,2}\geq 0$. Then, after $u_{2k}$ is determined, we set 
\[ A_{2k+1} = \begin{bmatrix}
1 & 0\\
\epsilon & 2
\end{bmatrix}   \]
with $\epsilon = 0.5 \frac{ |x_{2k,2}| }{\max(|x_{2k,1} |,1)}\mathrm{sign}(u_{2k})$
 Then, after $(A_{2k+1}, B_{2k+1})$ is revealed to the learner, it takes an action $u_{2k+1}$, resulting in $x_{2k+2} = A_{2k+1} A_{2k} x_{2k} + A_{2k+1} B_{2k} u_{2k} + B_{2k+1} u_{2k+1}$. Since $A_{2k} = I$ and $B_{2k} = B_{2k+1} = [1,0]^\top$, the second coordinate of $x_{2k+2}$ can then be written as
 \[x_{2k+2,2} = \epsilon x_{2k,1} + 2 x_{2k,2} + \epsilon u_{2k}. \]
 By the way $\epsilon$ is chosen, we have $\epsilon u_{2k} \geq 0$, and $|\epsilon x_{2k,1}|\leq 0.5 | x_{2k,2}| $, and therefore, $x_{2k+2,2} \geq 1.5 x_{2k,2}$. 
 
 Since the system starts with $x_0 = [1,1]^\top$, applying the argument above recursively, we have for any $k$, $x_{2k,2} \geq 1.5^k $, i.e. the state blows up. 
 \end{example}

\section{Proof of Theorem~\ref{thm:predictions}}\label{sec:predictions}
The feasibility directly follows Lemma~\ref{lemma:feasibility}. For stability, recall that we have rewritten the dynamics as,
\begin{align*}
x_{t+H} = \tilde{A}_t x_t + \tilde{B}_t \begin{bmatrix}u_{t+H-1}\\ u_{t+H-2}\\ \vdots \\u_{t}\end{bmatrix}  +  \underbrace{[I, A_{t+H-1}, \cdots,  A_{t+H-1}...A_{t+1}]\begin{bmatrix}w_{t+H-1}\\ w_{t+H-2}\\ \vdots \\w_{t}\end{bmatrix} }_{:= \tilde{w}_t}
\end{align*}
where we have defined $\tilde{A_t}:= A_{t+H-1}\cdots A_t$, $\tilde{B_t} := [B_{t+H-1}, A_{t+H-1}B_{t+H-2}, \cdots,  A_{t+H-1}...A_{t+1} B_t]$. By Theorem~\ref{thm:stability}, we have, $\forall k\geq 1$
\begin{align*}
    \Vert x_{kH+1}\Vert \leq \rho^{k}\Vert x_{1}\Vert + \frac{\kappa\rho}{1-\rho} \sup_{0\leq \tau < k} \Vert \tilde{w}_{\tau H+1}\Vert. 
\end{align*}
Note that, for any $\tau$,
\begin{align*}
    \Vert \tilde{w}_{\tau H+1}\Vert &\leq \sum_{\ell =1}^H  \Vert A_{\tau H+H} \Vert \cdots \Vert A_{\tau H+\ell+1}\Vert \Vert w_{\tau H + \ell }\Vert\\
    &\leq (1+a+\ldots + a^{H-1}) \sup_{ \tau H+1 \leq s \leq (\tau+1)H} \Vert w_s\Vert:= \kappa_A \sup_{ \tau H+1 \leq s \leq (\tau+1)H} \Vert w_s\Vert.
\end{align*}
These show that, 
\begin{align*}
    \Vert x_{kH+1}\Vert \leq \rho^{k}\Vert x_{1}\Vert + \frac{\kappa_A \kappa\rho}{1-\rho} \sup_{ 1 \leq s \leq  k H} \Vert w_{s}\Vert. 
\end{align*}
The above already shows the boundedness of states for time indexes $t=1\mod(H)$. To show the boundedness of the states for all $t$, we introduce the following lemma, whose proof is deferred to the end of this section. 
\begin{lemma}\label{lem:prediction:connecttoendpoint}
For any $k$, and $1\leq \ell\leq H$, we have, 
\begin{align*}
    \Vert x_{kH+ \ell} \Vert \leq \kappa_A' \Vert x_{kH+1} \Vert + \kappa_A \sup_{kH+1\leq s<kH+\ell} \Vert w_s\Vert. 
\end{align*}
where $\kappa_A' = a^{H-1} + b^2(1+a+\cdots + a^{H-1})^2 \kappa_R \frac{ \kappa + a^{H}}{\sigma}$. 
\end{lemma}
For any $t$, let $k$ be such that $t = kH+\ell$ with $1\leq \ell \leq H$. Then, using Lemma~\ref{lem:prediction:connecttoendpoint}, we have, 
\begin{align*}
    \Vert x_{t}\Vert &\leq \kappa_A' \Vert x_{kH+1} \Vert + \kappa_A \sup_{kH+1\leq s < t}\Vert w_s\Vert \\
    &\leq \kappa_A'  \rho^{k}\Vert x_{1}\Vert +\kappa_A' \frac{\kappa_A \kappa\rho}{1-\rho} \sup_{1 \leq s \leq  kH} \Vert w_s\Vert + \kappa_A \sup_{kH+1\leq s < t}\Vert w_s\Vert \\
    &\leq \kappa_A'  \rho^{\frac{t}{H} -1 }\Vert x_{1}\Vert +\kappa_A'\kappa_A \kappa \max(1, \frac{\rho}{1-\rho}) \sup_{1 \leq s \leq  t} \Vert w_s\Vert  . 
\end{align*}
So Theorem~\ref{thm:predictions} is proven. In the rest of the section, we provide a proof for Lemma~\ref{lem:prediction:connecttoendpoint}. 

\bigskip
\noindent\textbf{Proof of Lemma~\ref{lem:prediction:connecttoendpoint}:}  Note that, by Lemma~\ref{lemma:seq_stability}, we have the policies $K_{kH+1}$ is $(\kappa,\gamma,\rho)$-sequential strongly stable, with $\gamma = 1-\sqrt{\alpha}$ and $\kappa = \frac{\kappa_W}{\sqrt{1-\alpha}} $. This implies
$\Vert \tilde{A}_{kH+1} + \tilde{B}_{kH+1} K_{kH+1}\Vert \leq (1-\gamma) \kappa  $, showing
$$ \Vert \tilde{B}_{kH+1} K_{kH+1}\Vert \leq (1-\gamma)\kappa + a^H, $$
which further leads to
\begin{align*}
\Vert K_{kH+1}\Vert \leq K_{\max} = \kappa_R \frac{b(1+a+\ldots+a^{H-1}) }{\sigma} ((1-\gamma)\kappa + a^H). 
\end{align*}
As such, 
\begin{align*}
    \Vert x_{kH+ \ell} \Vert &\leq \Vert A_{kH+\ell-1} \cdots A_{kH+1}\Vert \Vert x_{kH+1}\Vert \\
    &+ \Vert \begin{bmatrix}
    B_{kH+\ell-1} , A_{kH+\ell-1} B_{kH+\ell-2},\ldots , A_{kH+\ell-1} \cdots A_{kH+2} B_{kH+1}
    \end{bmatrix} \begin{bmatrix} u_{kH+\ell-1}\\\vdots \\ u_{kH+1}\end{bmatrix}\Vert\\
    &+ \Vert [I, A_{t+H-1}, \cdots,  A_{kH +\ell-1 }...A_{kH+1}] \begin{bmatrix} w_{kH+\ell-1}\\\vdots \\ w_{kH+1}\end{bmatrix}\Vert\\
    &\leq a^{\ell-1} \Vert x_{kH+1}\Vert + b(1+a+\cdots+a^{\ell-2}) \Vert K_{kH+1}\Vert \Vert x_{kH+1}\Vert+ (1+a+ \cdots + a^{\ell-1}) \sup_{kH+1\leq s< kH+\ell} \Vert w_s\Vert\\
    &< \kappa_A' \Vert x_{kH+1} \Vert + \kappa_A \sup_{kH+1\leq s<kH+\ell} \Vert w_s\Vert. 
\end{align*}
\qed

\section{Proof of Theorem~\ref{thm:estimation}}
\label{sec:estimation}
\begin{proof} The proof is divided by two parts. For the first part, we prove the ISS property given the upper bound $K_{max}$ on the controllers $\Vert K_t\Vert$. In the second part, we provide such an upper bound $K_{max}$.

\noindent\textbf{Proof of ISS. }By Lemma~\ref{lem:iss}, to show ISS we only need to show that $(K_t)_{t=0}^\infty$ is sequential strongly stable for system $(A_t,B_t)_{t=0}^\infty$. By assumption, $\{ K_t\}_{t=0}^{\infty}$ is $(\kappa, \gamma, \rho)$-sequential strongly stable for the system $(\hat{A}_t,\hat{B}_t)_{t=0}^\infty$ with $(\kappa, \gamma, \rho)$ defined by Lemma~\ref{lemma:seq_stability} as $\kappa=\frac{\kappa_W}{\sqrt{1-\alpha}},\gamma=1-\sqrt{\alpha}, \rho=\sqrt{\frac{\alpha}{1-\alpha}}$ where $\kappa_W = \Vert W\Vert \Vert W^{-1}\Vert$.
Thus, there exist matrices $H_1, H_2, ..., $ and $L_1,L_2 ..., $ such that $\hat{M}_t:= \hat{A}_t+ \hat{B}_t K_t = H_t L_t H_t^{-1}$ with the following properties: (a) $\Vert L_t\Vert \leq 1-\gamma$; (b) $\Vert H_t\Vert \leq \beta_1$ and $\Vert H_t^{-1}\Vert  \leq 1/\beta_2$ with $\kappa = \beta_1/{\beta_2}$; (c)  $\Vert H_{t+1}^{-1} H_t\Vert \leq \frac{\rho}{1-\gamma}$. 
 With this decomposition for $\hat{M}_t$, we show that $M_t:= A_t + B_tK_t$ can be decomposed similarly. Let $\Delta_t = M_t - \hat{M}_t$. Then, 
\begin{equation}
    M_t = \hat{M}_t+\Delta_t = H_t L_t H_t^{-1} + H_tH_t^{-1} \Delta_t H_t H_t^{-1} = H_t (L_t+H_t^{-1} \Delta_t H_t) H_t^{-1} = H_t L_t'H_t^{-1}. \label{eq:estimation_M_decomp}
\end{equation}
We next upper bound $\Vert L_t'\Vert$. Notice that
\begin{align*}
    \Vert\Delta_t\Vert &= \Vert M_t-\hat{M}_t\Vert  = \Vert A_t + B_t K_t - \hat{A}_t - \hat{B}_t K_t\Vert \\
    & \leq \Vert A_t -\hat{A}_t\Vert  + \Vert B_t-\hat{B}_t\Vert  \Vert K_t\Vert  \\
    & \leq \max\{\Vert A_t -\hat{A}_t\Vert , \Vert B_t-\hat{B}_t\Vert \}(1+\Vert K_t\Vert ) \\
    & \leq \delta \frac{\gamma}{\kappa(1+K_{max})} (1+\Vert K_t\Vert )  \\
    & \leq \delta \frac{\gamma}{\kappa}.
\end{align*}
Then, we have the following bound on $\Vert L_t'\Vert$, 
\begin{align}
    \Vert L_t'\Vert  &= \Vert L_t+H_t^{-1} \Delta_t H_t\Vert  \nonumber\\
    & \leq \Vert L_t\Vert  + \Vert H_t^{-1}\Vert \Vert \Delta_t\Vert \Vert H_t\Vert  \nonumber\\
    & \leq 1-\gamma+ \delta \gamma = 1-(1-\delta)\gamma. \label{eq:estimation_Lprime}
\end{align}
Define $\gamma' = (1-\delta)\gamma$ and $\rho' =\frac{1-(1-\delta) \gamma}{1-\gamma} \rho$. 
We claim that $(K_t)_{t=0}^\infty$ is $(\kappa, \gamma', \rho')$ sequential strongly stable for system $(A_t,B_t)_{t=0}^\infty$. In the decomposition~\eqref{eq:estimation_M_decomp}, $L_t'$ satisfies $\Vert L_t'\Vert\leq 1-\gamma'$ (which we have proved in \eqref{eq:estimation_Lprime}), and by definition $ H_t $ satisfies $\Vert H_t\Vert \leq \beta_1$ and $\Vert H_t^{-1}\Vert  \leq 1/\beta_2$ with $\kappa = \beta_1/{\beta_2}$;  $\Vert H_{t+1}^{-1} H_t\Vert \leq \frac{\rho}{1-\gamma}= \frac{\rho'}{1-\gamma'}$. Therefore,
the only conditions we need to verify are (a) $0<\gamma' \leq 1$ and (b) $0 \leq \rho' < 1$. Condition (a) is trivial since for any $\alpha \in [0, \frac{1}{2})$, we have $\delta \in (0, \frac{\sqrt{1-\alpha}-\sqrt{\alpha}}{1-\sqrt{\alpha}})\subset (0, 1)$ and hence $\gamma' = (1-\delta)\gamma$ lies in $(0, 1]$. For condition (b), as $\rho' =\frac{1-(1-\delta) \gamma}{1-\gamma} \rho$, where $\rho = \sqrt{\frac{\alpha}{1-\alpha}}$ and $\gamma = 1-\sqrt{\alpha}$, it follows that $\rho'<1$ is equivalent to the following condition:
\begin{align}
    \rho' &= \frac{1-(1-\delta) \gamma}{1-\gamma} \rho =  \frac{\delta + (1-\delta) \sqrt{\alpha}}{\sqrt{1-\alpha}} < 1.
\end{align} 
Reorganizing the terms we have, 
\begin{align*}
    \delta< \frac{\sqrt{1-\alpha}-\sqrt{\alpha}}{1-\sqrt{\alpha}},
\end{align*}
which is satisfied by our condition on $\delta$.


 \noindent\textbf{Upper bounding $\Vert K_t\Vert$.} The only thing remains to show is that there exists an upper bound on the controller gain matrix under the conditions stated in the Theorem. 


Note that $K_t$ must satisfy $\hat{M}_t = \hat{A}_t + \hat{B}_t K_t$. On the other hand,
by the COCO-LQ formulation, $K_t$ solves the following problem
\begin{subequations}
  \begin{align*}
    \min_{K_t} &\qquad \Tr(R K_t \Sigma_{xx} K_t^\top ) \\
    \text{s.t. } & \hat{M}_t = \hat{A}_t + \hat{B}_t K_t.
\end{align*}
\end{subequations}
where $\Sigma_{xx}$ is the solution to the COCO-LQ formulation at time $t$. The Lagrangian of the above optimization is
\begin{align*}
    L(K_t, \Gamma) =\Tr(R  K_t \Sigma_{xx} K_t^\top) + \Tr(\Gamma^\top  (\hat{A}_t + \hat{B}_t K_t- \hat{M}_t)).
\end{align*}
The optimizer $K_t$ must satisfy the following condition,
$$\nabla_{K_t} L(K_t, \Gamma) = 2 R K_t \Sigma_{xx}   + \hat{B}_t^\top  \Gamma = 0.$$
$$\hat{M}_t = \hat{A}_t + \hat{B}_t K_t.$$
Together, we obtain $K_t = R^{-1} \hat{B}_t^\top  (\hat{B}_t R^{-1} \hat{B}_t^\top )^{-1} (\hat{M}_t - \hat{A}_t)$ and
\begin{align}\label{eq:case3}
    \Vert K_t\Vert  \leq  \Vert R^{-1}\Vert  \Vert \hat{B}_t^\top \Vert  \frac{1}{ \sigma_{\min}( R^{-1})   \sigma_{\min}(\hat{B}_t \hat{B}_t^\top  ) } (\kappa(1-\gamma) + \Vert \hat{A}_t\Vert ),
\end{align}
where we have used by the $(\kappa,\gamma,\rho)$ sequential strong stability of $(K_t)_{t=0}^\infty$ for system $(\hat{A}_t,\hat{B}_t)_{t=0}^\infty$, we have $\Vert \hat{M}_t\Vert = \Vert \hat{A}_t + \hat{B}_t K_t\Vert \leq \kappa (1-\gamma)$. By our condition, we have $\Vert \hat{A}_t\Vert  \leq \bar{\sigma}_A$, $  \Vert \hat{B}_t\Vert  \leq \bar{\sigma}_B$, $\hat{B}_t \hat{B}_t^\top \succeq \underline{\sigma}_B^2$, and $\kappa_R = \Vert R\Vert \Vert R^{-1}\Vert$.  Then, we have $\Vert K_t\Vert $ upper bounded as follows,
$$\Vert K_t\Vert  \leq \kappa_R \frac{\bar{\sigma}_B}{\underline{\sigma}_B^2}(\kappa (1-\gamma) + \bar{\sigma}_A).$$ 
\end{proof}

\revise{
\section{Supplementary Experiment Results: Nonlinear System Control via Local Linear Approximation}
\label{sec:nonlinear}
We further test the performance of COCO-LQ for online control of linear time varying systems that are derived from local linearization of nonlinear systems. In particular, we consider the pendulum swingup task, where the system dynamics are described by,
\begin{equation}
\label{eq:pendulum}
    \begin{bmatrix}
    \dot{\theta} \\
    \ddot{\theta}
    \end{bmatrix} = \begin{bmatrix}
      \dot{\theta} \\
      \frac{g}{l} sin \theta 
    \end{bmatrix} + 
    \begin{bmatrix}
    0 \\
    \frac{1}{ml^2}
    \end{bmatrix} u = f(x, u)
\end{equation}
where $x = [\theta, \dot{\theta}]$ represents the system state, in which $\theta$, $\dot{\theta}$ and $\ddot{\theta}$ are the angle, angular velocity and angular acceleration.  $g$ is the gravitational acceleration, $l$ and $m$ are the length and mass of the pendulum. The control goal is to stabilize the pendulum at the straight up position with $\theta = 0$ and $\dot{\theta} = 0$, starting from any initial angle and velocity. The pendulum dynamics described by Eq~\eqref{eq:pendulum} is a nonlinear system, by taking local linear approximation at each $x_t$, we can approximate the nonlinear system via a linear time varying system,
\begin{equation}
    \begin{bmatrix}
    \dot{\theta} \\
    \ddot{\theta}
    \end{bmatrix} = \underbrace{\begin{bmatrix}
    0 & 1\\
    \frac{g}{l} \cos \theta_t & 0
    \end{bmatrix}}_{A_t}
    \begin{bmatrix}
    \theta \\
    \dot{\theta}
    \end{bmatrix} + 
    \underbrace{\begin{bmatrix}
    0 \\
    \frac{1}{ml^2}
    \end{bmatrix}}_{B_t} u 
\end{equation}

\begin{figure}[htbp]
\includegraphics[width=0.24\columnwidth]{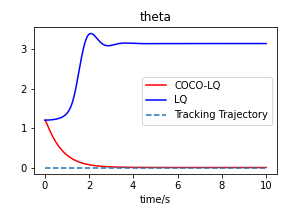}
\includegraphics[width=0.24\columnwidth]{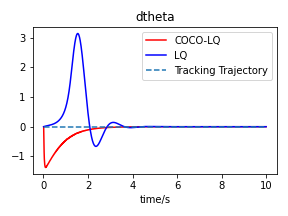}
\includegraphics[width=0.24\columnwidth]{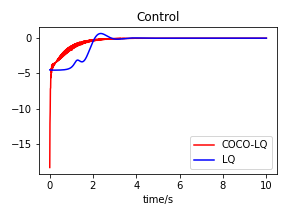}
\includegraphics[width=0.24\columnwidth]{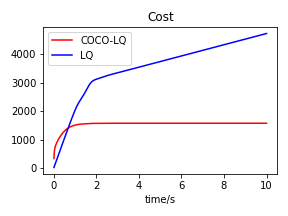}
\caption{Performance comparison of COCO-LQ and LQ on inverted pendulum via locally linearization. The left two figures show the state evolution of angle $\theta$ and angular velocity $\dot{\theta}$. Initial angle is set as $\theta = 1.2 \rad$, and the desired state is $\theta = \dot{\theta} = 0$. The right two figures show the control action and cost comparison, with $Q = R = I$.}
\label{fig:pendulum_result}
\end{figure}

Fig~\ref{fig:pendulum_result} compares the performance of the proposed COCO-LQ approach and baseline LQ approach. As we can observed from the left two plots, COCO-LQ is able to stabilize the pendulum at the desired position due to its robustness under model estimation error (caused by local linear approximation) and time-varying dynamics, while the naive LQR approach fails to achieve the swingup task and stabilize the system. The right two plots show the evolution of the control efforts and cost. COCO-LQ initially outputs significantly larger control actions compared to the baseline LQ controller, to stop the pendulum from falling over. The overall control cost of COCO-LQ quickly converges once the swingup task is finished, while the cost of LQ control keeps growing.}

\end{document}